\documentclass[paper]{amsart}

\usepackage{amsfonts}
\usepackage{amssymb}
\usepackage{graphicx}

\def\f{\varphi}               \def\z{\zeta}

\def\D{{\mathbb D}}  
  
\def\R{{\mathbb R}}

\def\({\left(}       \def\){\right)}

\newtheorem{prop}{\sc Proposition}
\newtheorem{lem}{\sc Lemma}
\newtheorem{thm}{\sc Theorem}

\newtheorem{other}{\sc Theorem}              
\newenvironment{pf}{\noindent{\textit{Proof. }}}{$\Box$ }

\begin{document}
\title[Generalized harmonic Koebe functions]
{Generalized harmonic Koebe functions}

\author[A. Ferrada-Salas]{\'Alvaro Ferrada-Salas}
\address{Facultad de Matem\'aticas, Pontificia Universidad Cat\'olica de Chile, Casilla 306, Santiago, Chile.} \email{alferrada@mat.puc.cl}

\author[M. J. Mart\'{\i}n]{Mar\'{\i}a J. Mart\'{\i}n}
\address{Department of Physics and Mathematics, University of Eastern Finland, P.O. Box 111, FI-80101 Joensuu, Finland.} \email{maria.martin@uef.fi}

\date{\today}
\thanks{The first author is supported by Conicyt grant 21120474, Chile. The second author is supported by
Academy of Finland grant 268009 and by Spanish MINECO Research Project MTM2012-37436-C02-02.}

\begin{abstract}

We present a family of sense-preserving harmonic mappings in the unit disk related to the classical generalized (analytic) Koebe functions. We prove that these are precisely the mappings that maximize simultaneously the real part of every Taylor coefficient as well as the growth and distortion of functions in affine and linear invariant families of complex-valued harmonic functions.

\end{abstract}
\maketitle


\section*{Introduction}\label{sec-introduction}

 Let $S$ be the family of all univalent (one-to-one) analytic mappings $\varphi$ in the unit disk $\D$ with the normalizations $\varphi(0)=1-\varphi'(0)=0$.
\par
It was as early as 1916 when Bieberbach \cite{Bieberbach1, Bieberbach2} proved that the second Taylor coefficient of any function in the class $S$ is bounded by $2$ and conjectured that the bound for the $n$-th Taylor coefficient of mappings in $S$ should be $n$. The Bieberbach conjecture was proved by De Branges \cite{DeBranges} in 1985.  The Koebe mapping
\begin{equation}\label{eq-koebe}
k(z)=\frac{z}{(1-z)^2}\,,\quad z\in\D\,,
\end{equation}
belongs to $S$ and has Taylor series expansion
\[
k(z)=z+\sum_{n=2}^\infty n z^n\,.
\]
Hence, the bound $n$ for the $n$-th Taylor coefficient of functions in $S$ is sharp.
\par
Notice that since whenever $\lambda\in\partial\D$ and $\varphi\in S$, the \emph{rotation} $\varphi_\lambda$ defined by $\varphi_\lambda(z)=\overline\lambda \varphi(\lambda z)\in S$, the problem of maximizing the modulus of the $n$-th Taylor coefficient of functions in $S$ is equivalent to the problem of maximizing the real part of such coefficient. In the latter case (as it is explained in \cite[Sec. 2.9]{Dur-Univ}), by the use of a variational technique, Marty \cite{Marty} proved that the relation
\begin{equation}\label{eq-Marty}
(n+1)a_{n+1}=2a_2a_n+(n-1)\ \overline{a_{n-1}}
\end{equation}
must be satisfied by the coefficients of each mapping $\varphi\in S$ whose $n-$th Taylor coefficient has maximum real part. By denoting the $n$-th Taylor coefficient of a function $\varphi\in S$ by $a_n(\varphi)$, we see that the Koebe function (\ref{eq-koebe}) has the property that
\[
\sup_{\varphi\in S} Re\{a_n(\varphi)\}=Re\{a_n(k)\}=a_n(k)
\]
for all $n\geq 2$. Obviously, the Koebe mapping satisfies Marty's relation (\ref{eq-Marty}) for all $n$.
\par\smallskip
The variational method used by Marty to obtain (\ref{eq-Marty}) works more generally in the setting of \emph{linear invariant families}; this is,  families of locally univalent analytic functions $\varphi$ in the unit disk normalized as above and which are closed under the transformation
$$K_{\zeta}(z)=\frac{\f\left(\displaystyle\frac{\zeta+z}{1+\overline{\zeta}z}\right)
-\f(\zeta)}{(1-|\zeta|^2)\f^\prime(\zeta)}\,,\quad \zeta\in\D.$$
\par
Several important properties, such as growth, covering, and distortion are determined by the \emph{order} of a linear invariant family $F$ defined by
$$\alpha(F)=\sup_{\varphi\in F} |a_2(\varphi)|=\sup_{\varphi\in F} Re\{a_2(\varphi)\}\,.$$
\par
\par
We refer the reader to the works \cite{POM-I, POM-II}, where Pommerenke studies and carries through a detailed analysis of linear invariant
families. We would also like to mention \cite[Ch. 5]{GK} as a good reference on the topic.
\par
Explicit examples of linear invariant families are the class $S$ defined above (see the books \cite{Dur-Univ} or \cite{P} for more details related to this class) and also the family $F_M$ of all normalized locally univalent analytic functions in the unit disk satisfying
\[
\sup_{|z|<1} |S\varphi(z)| (1-|z|^2)^2\leq M\,,
\]
where $S\varphi$ is the Schwarzian derivative of $\varphi$:
\[
S\varphi=\left(\frac{\varphi''}{\varphi'}\right)'-\frac 12 \left(\frac{\varphi''}{\varphi'}\right)^2\,.
\]
\par
Regarding these families $F_M$, Pommerenke \cite{POM-I} proves that
$$\sup_{\varphi\in F_M} |a_2(\varphi)|=\sqrt{1+\frac{M}{2}}\,.$$
It is a straightforward calculation to show that given $M\geq 0$, the function
\begin{equation}\label{eq-koebegeneral}
k_a(z)=\frac{1}{2a}\left[\left(\frac{1+z}{1-z}\right)^a-1\right]\,, \quad |z| < 1\,,
\end{equation}
where $a=\sqrt{1+M/2}\,,$ belongs to $F_M$ and satisfies
\[
\frac 12 |k_a''(0)| = \sqrt{1+\frac{M}{2}}\,.
\]
\par\smallskip
These functions $k_a$ are called \emph{generalized Koebe functions} and they are known to be extremal for a number of problems (see, for instance, \cite[Ch.~ 11]{Goodman} or \cite{Koepf}). In particular, in the next proposition we point out explicitly one property satisfied by $k_a$. Though the proof is not difficult and the result is probably known to the experts, we have not been able to find an explicit reference. We include it here for the sake of completeness.
\begin{prop}\label{prop-analytic}
Let $F$ be a linear invariant family of analytic functions in the unit disk. Assume that there exists a function $\Phi \in F$ such that for all positive integer $n,\ n\geq 2$,
\[
\sup_{\varphi\in F} Re \{a_n(\varphi)\}=Re\{a_n(\Phi)\}\,.
\]
Then, $\Phi=k_a$, where $a$ is the order of $F$.
\end{prop}
\begin{pf}
In order to simplify the notation, let us use $A_n$ to denote the Taylor coefficients of $\Phi$; this is $A_n=a_n(\Phi)$. We first note that all the coefficients $A_n$ must be non-negative real numbers since, otherwise, we could use a rotation of the form $\Phi_\lambda(z)=\overline \lambda \Phi(\lambda z)$ (with appropriate $\lambda\in\partial \D$ that could depend on $n$) to obtain another function in $F$ with bigger Taylor coefficient than that of $\Phi$. In fact, the coefficient $A_2$ is strictly positive since it coincides with the order of $F$ and, as Pommerenke showed in \cite{POM-I}, the order of any linear invariant family of analytic functions is always greater than or equal to $1$.
\par
Now, using the assumption that $\Phi$ maximizes the real part of any Taylor coefficient, we get that the Marty relation (\ref{eq-Marty}) holds for all $n\geq 1$. Hence, using also that $A_n \in\R$, we have that for all such $n$,
\[
(n+1)A_{n+1}=2A_2A_n+(n-1)A_{n-1}\,.
\]
Therefore, we obtain
\begin{eqnarray*}
	\Phi'(z)&=&1+\sum_{n=2}^{\infty}nA_nz^{n-1}\\
	&=&1+\sum_{n=1}^{\infty}(n+1)A_{n+1}z^{n}\\
	&=&1+2A_2\sum_{n=1}^{\infty}A_nz^n+\sum_{n=1}^{\infty}(n-1)A_{n-1}z^{n}\\
	&=&1+2A_2 \Phi(z)+z^2 \Phi'(z)\,.
\end{eqnarray*}
Thus, $\Phi$ solves the linear differential equation
\begin{equation}\label{eq-edo}
(1-z^2)\phi'(z)=1+\alpha\phi(z)\,,\quad z\in\D\,,
\end{equation}
with $\alpha=2A_2$. Using basic techniques for solving linear differential equations of first order, it is easy to see that the (unique) solution to (\ref{eq-edo}) with initial datum $\phi(0)=0$ is $\phi=k_{a}$ with $a=\alpha/2=A_2$. This ends the proof.
\end{pf}
\par\smallskip
Notice that Proposition~\ref{prop-analytic} can be re-stated as follows:

\begin{quote}
\emph{The unique possible function maximizing the real part of every Taylor coefficient of functions in a linear invariant family of analytic mappings in $\D$ is a generalized Koebe function of the form (\ref{eq-koebegeneral})}.
\end{quote}
\par\medskip
In this paper, we will determine what are the \emph{harmonic} functions maximizing the Taylor coefficients in affine and linear invariant families of harmonic mappings in the unit disk. In other words, we will obtain an analogue of Proposition~\ref{prop-analytic} in the harmonic setting. We will see that those functions (we call them \emph{generalized harmonic Koebe functions}) for which the real part of any Taylor coefficient is maximum among all the mappings in affine and linear invariant families can be defined via the method of `shear construction' developed by Clunie and Sheil-Small \cite{CSS}. The analytic generalized Koebe functions of the form (\ref{eq-koebegeneral}) will play an important role in the construction too.

\section{Background}\label{sec-preliminar}
Note that since the unit disk $\D$ is simply connected, any complex-valued harmonic mapping $f$ in $\D$ has a canonical decomposition $f=h+\overline{g}$, where $h$ and $g$ are analytic in $\D$. Following \cite{CSS}, we call $h$ the \emph{analytic part of $f$} and $g$ the \emph{co-analytic part of $f$}. The harmonic mapping $f$ is analytic if and only if $g$ is constant. A harmonic mapping $f=h+\overline g$ is \emph{sense-preserving} if it has positive Jacobian; this is if $h^\prime$ does not vanish in the unit disk and the (second complex) \emph{dilatation} $\omega=g'/h'$ has the property that $|\omega|<1$ in $\D$. We refer the reader to the book by Duren \cite{Dur-Harm} for an excellent exposition on harmonic mappings.
\par
\subsection{Affine and linear invariant families of harmonic mappings} \label{ssec-alfamilies}
Let $F_H$ be a family of sense-preserving harmonic mappings $f=h+\overline{g}$ in $\D$, normalized by $h(0)=g(0)=0$ and $h'(0)=1$. The family is said to be {\it affine and linear invariant} if it
closed under \emph{Koebe transforms} 
\begin{equation}\label{eq-koebetransform}
K_{\zeta}(f)(z)=\frac{f\left(\displaystyle\frac{z+\zeta}{1+\overline{\zeta} z}\right)-f(\zeta)}{(1-|\zeta|^2)h'(\zeta)}\,,\quad |\zeta|<1\,,
\end{equation}
and \emph{affine changes}
\begin{equation}\label{eq-affinetransform}
A_\varepsilon(f)(z)=\frac{f(z)-\overline{\varepsilon f(z)}}{1-\overline\varepsilon g'(0)}\,,\quad |\varepsilon|<1\,.
\end{equation}
\par
We remark that we are assuming that all the members in an affine and linear invariant mapping of harmonic functions are sense-preserving harmonic mappings in the unit disk.
\par
For a given affine and linear invariant family of harmonic mappings $F_H$, we use $F_H^0$ to denote the subset of functions $f=h+\overline g\in F_H$ satisfying $g'(0)=0$. Note that $F_H^0$ is not an affine and linear invariant family in general.
\par
Sheil-Small \cite{S-S} offers a deep study of affine and linear invariant families $F_H$ of harmonic mappings in $\D$. The \emph{order} of the affine and linear invariant family, given by
$$ \alpha(F_H)=\sup_{f=h+\overline g\in F_H}|a_2(h)|=\frac 12 \sup_{f=h+\overline g\in F_H}|h''(0)| \,,$$
plays once more a special role in the analysis.
\par\smallskip
A  special example of affine and linear invariant family is the class $S_H$ of (normalized) sense-preserving harmonic mappings which are univalent in the unit disk (see \cite{CSS} and \cite[Ch. 5]{Dur-Harm}). Also, in \cite{HM2}, the authors introduce a definition for the Schwarzian derivative $S_f$ for locally univalent harmonic mappings. Using this definition for the Schwarzian derivative, it is proved in \cite{CHM} that the family $F_H^M$ of sense-preserving harmonic mappings $f=h+\overline{g}$ in $\D$, with $h(0)=g(0)=0,\ h'(0)=1$ and $||S_f||\leq M$, is affine and linear invariant.
\subsection{Harmonic Marty relations}\label{ssec-Marty}
As mentioned in \cite[p. 101]{Dur-Harm}, a slight modification of the Marty variation in the analytic case leads to analogues of the Marty relation (\ref{eq-Marty}) for harmonic mappings. We refer the reader to \cite[Sec. 6.5]{Dur-Harm} for a proof of the fact that whenever $F_H$ is an affine and linear invariant family of harmonic mappings and there is a function $f=h+\overline g \in F_H$ with
\[
h(z)=z+\sum_{n=2}^\infty a_n z^n \quad {\rm and} \quad g(z)=\sum_{n=2}^\infty b_n z^n
\]
that maximizes the real part of the $n$-th Taylor coefficient of analytic parts of functions in $F_H^0$, then
\begin{equation}\label{eq-martyan}
(n+1)a_{n+1}=2a_2a_n+2b_2\overline{b_n}+(n-1)\overline{a_{n-1}}
\end{equation}
must hold. Similarly, if $f$ has a coefficient $b_n$ of maximum real part, then
\begin{equation}\label{eq-martybn}
(n+1)b_{n+1}=2a_2\overline{b_n}+2b_2a_n+(n-1)\overline{b_{n-1}}\,.
\end{equation}
\par\smallskip
Perhaps at this point we should stress that as long as a given function $f=h+\overline g\in F_H^0$ (with $h$ and $g$ as above) has a coefficient $a_n$ (\emph{resp.} $b_n$) of maximum real part, then indeed, $a_n$ (\emph{resp.} $b_n$) must be a non-negative real number since, otherwise, we can consider the function $f_\lambda$ defined by
\[
f_\lambda(z)=\overline\lambda f(\lambda z) = \overline\lambda h(\lambda z)+\overline{\lambda g(\lambda z)}\,,
\]
with an appropriate $\lambda,\ |\lambda| =1$, to get another function in $F_H^0$ with bigger coefficient(s) than those of $f$.

\subsection{The shear construction}\label{ssec-shear}
An effective and beautiful way of constructing sense-preserving univalent harmonic mappings in the unit disk is the so called \emph{shear construction} introduced in the paper \cite{CSS} by Clunie and Sheil-Small. It is based on the following result. Recall that a function $f$ is convex in the $\theta$ direction ($0  \leq \theta < \pi$) if the intersection of the domain $f(\mathbb{D})$ with every line parallel to the line through $0$ and $e^{i\theta}$  is connected.

\begin{other}\label{thm-shear}
Let $f=h+\overline g$ be a locally univalent harmonic mapping in the unit disk. Then it is univalent and convex in the $\theta$ direction if and only if the analytic function $h-e^{2i\theta} g$ is univalent and convex in the $\theta$ direction.
\end{other}
\par
The method of constructing univalent harmonic mappings using the shear construction can be, then, stated as follows. Consider any analytic function $\varphi$ in $\D$ which maps the unit disk onto a domain convex in the direction $\theta$. Also, let $\omega$ be any analytic function in the unit disk with $\omega(\D)\subset \D$. The \emph{harmonic shear in the $\theta$ direction of the function $\varphi$ with dilatation $\omega$} is the harmonic function $f=h+\overline g$, where $h$ and $g$ solve the linear system of equations
\begin{eqnarray}\label{eq-shear}
\left\{ \begin{array}{c} h-e^{2i\theta}g=\varphi  \\ g^\prime/h^\prime=\omega  \end{array}\right.\,,\quad h(0)=g(0)=0\,.
\end{eqnarray}
\par\smallskip
In \cite{CSS}, the authors gave two explicit examples of univalent harmonic mappings in $\mathbb{D}$ produced using this shear construction: the \emph{harmonic Koebe function} and the \emph{half-plane harmonic mapping}.

\par\smallskip
The harmonic Koebe function is defined by $K=H+\overline G$, where
\[
\quad H(z)=\frac{z-\frac 12 z^2+\frac 16 z^3}{(1-z)^3}\quad{\rm and}\quad G(z)=\frac{\frac 12 z^2+\frac 16 z^3}{(1-z)^3} \quad (z\in\D)\,.
\]
Note that $H$ and $G$ are the solution of the linear system of equations
\begin{eqnarray*}
\left\{ \begin{array}{c} H(z)-G(z)=k(z)  \\ G^\prime(z)/H^\prime(z)=z  \end{array}\right.\,, \quad z\in\mathbb{D}\,,
\end{eqnarray*}
with $H(0)=G(0)=0$. In other words, the harmonic Koebe function coincides with the horizontal shear (that is, the shear in the $0$ direction) of the classical analytic Koebe function with dilatation $\omega(z)=z$.
\par\smallskip
The second example constructed in \cite{CSS} using the shear construction is the half-plane harmonic mapping $L=h+\overline g$. It is the vertical shear (the shear in the $\pi/2$ direction) of the function $\ell (z)=z/(1-z)$ with dilatation $\omega(z)=-z$. That is, $h$ and $g$ solve the linear system
\begin{eqnarray*}
\left\{ \begin{array}{c} h(z)+g(z)=\ell(z)  \\ g^\prime(z)/h^\prime(z)=-z  \end{array}\right.\,, \quad z\in\mathbb{D}\,,
\end{eqnarray*}
$h(0)=g(0)=0$.
\par\smallskip

We finish this section by pointing out that according to Theorem~\ref{thm-shear}, the harmonic shear obtained by solving the system (\ref{eq-shear}) in the case when the function $\varphi$ is convex in the $\theta$ direction is univalent. It is not difficult to check that in the more general case when the function $\varphi$ is just supposed to be locally univalent, the solution $(h,g)$ to the linear system (\ref{eq-shear}) produces a sense-preserving (hence, locally univalent) harmonic function $f=h+\overline g$.
\section{Generalized harmonic Koebe functions} \label{ssec-harmkoebe}
Let $a$ be any complex number different from $0$ and consider the \emph{generalized Koebe function} $k_a$ as in (\ref{eq-koebegeneral}); this is,
\[
k_a(z)=\frac{1}{2a} \left[\left(\frac{1+z}{1-z}\right)^{a}-1\right]\quad (z\in\D)\,,
\]
where the branch of the logarithm is chosen so that $\log 1=0$ in
\[
\left(\frac{1+z}{1-z}\right)^{a}=\exp\left(a\log\frac{1+z}{1-z}\right)\,.
\]
Note that since
\[
k_a'(z)=\left(\frac{1+z}{1-z}\right)^{a}\cdot \frac{1}{1-z^2}\,,
\]
all these functions $k_a$ are locally univalent in the unit disk. Moreover, the limit (in compact subsets of $\D$) of the functions $k_a$ as $a\to 0$ coincides with
\begin{equation}\label{eq-k0}
k_0(z):=\frac 12\log\frac{1+z}{1-z}\,,
\end{equation}
a function that turns out to be univalent in $\D$.
\par
Hille, using generalized Koebe functions with $a=i\varepsilon$, where $\varepsilon$ is a real number with modulus small enough, proved that the constant $2$ in the well-known criterion of univalence in terms of the Schwarzian derivative due to Nehari is sharp (see \cite{Hille}). Indeed, in this paper \cite{Hille} Hille  proves that $k_a$ is univalent if and only if either $a$ or $-a$ belong to the (closed) disk centered at the point $z=1$ and with radius $1$. This implies, in particular, that for real values of $a$, $k_a$ is univalent if and only if $-2\leq a\leq 2$. (Note that $k_2$ equals the Koebe function (\ref{eq-koebe}).)
\par\medskip
In order to define an appropriate analogue of generalized Koebe functions in the harmonic setting, it seems to make sense to consider mappings $K=h+\overline g$ produced by the shear (in some direction $\theta$) of $k_a$. The following proposition somehow justifies this election. Recall that by Proposition~\ref{prop-analytic}, we know that if a function $\varphi$ \emph{in a linear invariant family of analytic functions} maximizes the real part of every Taylor coefficient of functions in the family, then $\varphi$ is a generalized Koebe function $k_a$ for some positive real number $a$. We should point out that, in principle, Proposition~\ref{prop-main} below  cannot be directly proved using Proposition~\ref{prop-analytic} since it is not clear that the family of analytic functions $\{h+\lambda g\colon h+\overline g \in F_H^0\}$ has to be, in general, a linear invariant family.

\begin{prop}\label{prop-main}
Let $F_H$ be an affine and linear invariant family of (normalized) sense-preserving harmonic mappings. Consider a complex number $\lambda$ with $|\lambda|=1$ and assume that there is a function $f_0=h_0+\overline{g_0}\in F_H^0$ such that
\begin{equation}\label{eq-propmain}
\sup_{f=h+g\in F_H^0} Re \{a_n(h)+\lambda b_n(g)\}=a_n(h_0)+\lambda b_n(g_0)
\end{equation}
for all $n\geq 2$. Then, $h_0+\lambda g_0$ equals the generalized analytic Koebe function $k_a$ as in (\ref{eq-koebegeneral}) with $a=a_2(h_0)+\lambda b_2(g_0)$.
\end{prop}
\smallskip\par
\begin{pf}
Using the normalizations of functions in $F_H$, it is obvious that (\ref{eq-propmain}) holds for $n=0$ and $n=1$. We use $A_n$ and $B_n$ to denote the coefficients $a_n(h_0)$ and $b_n(g_0)$, respectively. Note that our hypotheses shows that for all $n\geq 0$, $A_n+\lambda B_n$ is a non-negative real number.
\par
Take an arbitrary point $\zeta\in\D$ and consider the transformations (\ref{eq-koebetransform}) and (\ref{eq-affinetransform}) to produce the mappings
$$f_\zeta=A_{\omega_{\zeta}(0)}(K_\z(f_0))=h_\zeta^*+\overline{g_\zeta^*}$$
that belong to $F_H^0$ since $F_H$ is affine and linear invariant. Here we are using $\omega_\zeta$ to denote the dilatation of the function $K_\z(f_0)$.
\par\smallskip
As it is shown on \cite[p.102]{Dur-Harm}, the Taylor coefficients $a_n^\ast$ and $b_n^\ast$ of $h_\zeta^*$ and  $g_\zeta^*$, respectively, satisfy
\[
 a_n^\ast=A_n+[(n+1)A_{n+1}-2A_2A_n]\zeta-[2\overline{B_2}B_n+(n-1)A_{n-1}]\overline{\zeta}+o(|\zeta|)
\]
and
\[
b_n^\ast=B_n+[(n+1)B_{n+1}-2B_2A_n]\zeta-[2\overline{A_2}B_n+(n-1)B_{n-1}]\overline{\zeta}+o(|\zeta|)\,.
\]
Using the extremal property of $f_0$, we get
\[
Re\{a_n^\ast+\lambda b_n^\ast\} \leq Re\{A_n+\lambda B_n\}=A_n+\lambda B_n\,.
\]

By arguing as in \cite[p. 102]{Dur-Harm}, and using also that for all $n\geq 0$ we have $\overline{A_n+\lambda B_n}=A_n+\lambda B_n$ and that $|\lambda|=1$, we obtain that the following equality must hold for the coefficients $A_n$ and $B_n$, $n\geq 1$:
\begin{eqnarray}\label{eq-propcoeff}
\nonumber (n+1)(A_{n+1}+\lambda B_{n+1}) &\quad  -\quad \ (n-1)(A_{n-1}+\lambda B_{n-1})\\
& \quad = 2 (A_2+\lambda B_2) (A_n+\overline{\lambda B_n})\,.
\end{eqnarray}
In particular, we see that $z_n=A_n+\overline{\lambda B_n}$ is a real number for all $n\geq 1$. Since, by hypotheses, $w_n=A_n+\lambda B_n$ is a real number as well, we obtain (just subtracting both numbers) that $Im\{\lambda B_n\}=0$. This shows that $\lambda B_n$ and (hence $A_n$ too) are real numbers for all such values of $n$. Therefore, we can re-write (\ref{eq-propcoeff}) as
\begin{eqnarray}\label{eq-propcoeff2}
\nonumber (n+1)(A_{n+1}+\lambda B_{n+1}) & \quad  -\quad (n-1)(A_{n-1}+\lambda B_{n-1})\\
& \nonumber = 2 (A_2+\lambda B_2) (A_n+\lambda B_n)
\end{eqnarray}
and using the same arguments as in the proof of Proposition~\ref{prop-analytic} we see that $h+\lambda g=k_{A_2+\lambda B_2}$, as was to be shown.
\end{pf}
\par\medskip
Harmonic shears of generalized Koebe functions have already been considered in the literature. For instance, as was explained above, the harmonic Koebe function $K$ is the horizontal shear of the analytic Koebe mapping $k=k_2$ with dilatation $\omega(z)=z$. The half-plane harmonic mapping, $L$, equals the shear in the vertical direction of $k_1$ with dilatation $\omega(z)=-z$. Indeed, vertical shears of $k_1$ with dilatation $\omega(z)=-\eta z$ (where $\eta$ is a complex number of modulus $1$) coincide with the functions $F(z,\eta)$ defined on \cite[p.~158]{GN}, which turn out to be extreme points of an important family related to functions mapping the unit disk onto the half-plane $H=\{z\colon Re\{z\}>-1/2\}$. Harmonic shears of $k_0$ have been proved to have also interesting properties (see \cite[Sec. 2]{GN}). Moreover, every generalized analytic Koebe function $k_a$ equals the shear (in any direction) of $k_a$ itself with dilatation $\omega\equiv 0$.
\par\smallskip
It was proved in \cite{CHM} that the functions $f(=f_{a,R})=h+\overline g$
defined by
\begin{equation}\label{eq-harmonicgeneralizedKoebe}
\left\{
\begin{array}{c}
h-g=k_a\\
\omega=g^\prime/h^\prime=l_R
\end{array}\right.\,,\quad\quad h(0)=g(0)=0\,,
\end{equation}
where for $0<R<1$ the \emph{lens-map} $l_R$ is defined by
\[
l_R(z)=\frac{\left(\frac{1+z}{1-z}\right)^R-1}{\left(\frac{1+z}{1-z}\right)^R+1}\,,\quad |z|<1\,,
\]
are extremal mappings for the problem of maximizing the real part of the second Taylor coefficient of $h$ for $f=h+\overline g\in\mathcal (F^M_H)^0$, where $F_H^M$ is the family of normalized harmonic mappings with Schwarzian norm bounded by $M$ mentioned in Section~\ref{ssec-alfamilies}. As we shall see below, the election of the mapping $l_R$ as the dilatation in the system (\ref{eq-harmonicgeneralizedKoebe}) will produce interesting properties for the harmonic shear $f=h+\overline g$ obtained from the solution of this system. Indeed, we should also mention the close relation between the lens-maps and the generalized analytic Koebe functions. Concretely, as the reader may check, for any value of $R,\ 0<R<  1$,
\[
l_R=\frac{Rk_R}{1+Rk_R}\,.
\]
\par\medskip
To unify the notation, write $\ell_0\equiv 0$ and $\ell_1(z)=z$ in the unit disk. Motivated by all the facts mentioned (Proposition~\ref{prop-main} included), we define a \emph{generalized harmonic Koebe function} $K_H=K_H(\lambda, a, \mu, R)$, where $|\lambda|=|\mu|=1$, $0\leq R\leq 1$, and $a$ is any complex number, as the function $K_H=h+\overline g$, where $h$ and $g$ solve the system
\begin{equation}\label{eq-harmKoebe3}
\left\{
\begin{array}{c}
h-\lambda g= k_a\\
\omega=g^\prime/h^\prime=\mu l_R
\end{array}\right.\,,\quad\quad h(0)=g(0)=0\,.
\end{equation}
\par
It is easy to see that the solution $(h,g)$ of (\ref{eq-harmKoebe3}) exists for any values of $\lambda$, $\mu$, $R$, and $a$ as above and gives rise to a sense-preserving (hence locally univalent) harmonic mapping $K_H=h+\overline g$.
\par\medskip
In the remaining part of the paper, we will focus on analyzing the properties of these generalize harmonic Koebe functions in the cases when the parameters $\lambda$, $\mu$, and $a$ are real numbers. Note that in this case,  different symmetry properties come out. For instance, since for any real number $a$ and for any $0\leq R\leq 1$ we have $k_{-a}(z)=-k_a(-z)$ and $l_R(z)=-l_R(-z)$, it is easy to check that
\begin{itemize}
\item[(i)] $K_H(1,a, 1, R)=K_H(-1,a+R,1,R)$\,.
\item[(ii)] $K_H(-1,a, -1, R)=K_H(1,a+R,-1,R)$\,.
\item[(iii)] $K_H(\lambda,-a, \mu, R)(z)=-K_H(\lambda,a,-\mu,R)(-z)$\,.
\end{itemize}
These symmetries show, in particular, that we can restrict ourselves to the analysis of generalized harmonic Koebe functions of the form $K(1, a, 1, R)$, that will be denoted by $K_{a, R}$. In other words, a harmonic mapping $f=h+\overline g$ is a generalized harmonic Koebe function of the form $K_{a, R}$ (where $a$ is a real number and $0\leq R\leq 1$) if $h$ and $g$ solve the system (\ref{eq-harmonicgeneralizedKoebe}).

\par\medskip
In the next theorem, we characterize those (real) values of $a$ and $R$ for which $K_{a, R}$ is univalent.
\begin{thm}\label{thm-univalent}
Let $a$ be any real number and $0\leq R\leq 1$. Then $K_{a, R}$ is univalent if and only if $-2\leq a\leq 2$.
\end{thm}
\begin{pf}
Consider first the case when $a=0$. Then $K_{0, R}=h+\overline g$, where $h$ and $g$ solve
\[
\left\{
\begin{array}{c}
h-g= k_0\\
\omega=g^\prime/h^\prime=l_R
\end{array}\right.\,,\quad\quad h(0)=g(0)=0\,,
\]
with $k_0$ as in (\ref{eq-k0}). Note that this function $k_0$ maps the unit disk onto a convex domain (which is convex in the 0 direction). By Theorem~\ref{thm-shear} we conclude that $K_{0, R}$ is univalent for all values of $R$ in $[0, 1]$.
\par
Now, assume that $a>0$. If $a\leq 2$, then $k_a$ is convex in the horizontal direction. Hence, another application of Theorem~\ref{thm-shear} shows that for any value of $R\in [0,1]$, the function  $K_{a, R}$ is univalent in the unit disk if $0<a\leq 2$.
\par
Let us prove now that if $a>2$ then $K_{a, R}$ is not univalent. Since the functions $k_a$ and $l_R$ have real Taylor coefficients and $h'=k_a'/(1-l_R)$, we have that the Taylor coefficients of $h$ are real too. Moreover, $g$ has this property as well since $g=h-k_a$.
\par
Let us write $a=2+\varepsilon$ for an appropriate $\varepsilon >0$. The transformation $z\to (1+z)/(1-z)$ maps the unit disk onto the right half-plane $Re\{z\}>0$. Hence, there exists $z_1\in\D$ with
\[
\frac{1+z_1}{1-z_1}=e^{i\frac{\pi}{2+\varepsilon}}\,,
\]
so that $k_a(z_1)$ is a negative real number. (Note that $z_1$ cannot be real.) Set $z_2=\overline{z_1}$. Bearing in mind that $k_a$ has real Taylor coefficients, we see that for all $z\in\D$ the identity $k_a(\overline z)=\overline{k_a(z)}$ holds. In particular, since $k_a(z_1)$ is a real number we have that $k_a(z_2)=k_a(z_1)$. In other words, we have found two different points $z_1$ and $z_2$ in $\D$ with
\[
k_a(z_1)=h(z_1)-g(z_1)=h(z_2)-g(z_2)=k_a(z_2)\,,
\]
which obviously gives $h(z_1)+g(z_2)=h(z_2)+g(z_1)$.
\par
Therefore (using also that $g$ has real Taylor coefficients) we get
\begin{eqnarray*}
f(z_1)&=&h(z_1)+\overline{g(z_1)}=h(z_1)+\overline{g(\overline{z_2})}=h(z_1)+g(z_2)\\
&=& h(z_2)+g(z_1)=h(z_2)+g(\overline{z_2})=h(z_2)+\overline{g(z_2)}=f(z_2).
\end{eqnarray*}
This proves that $K_{a, R}$ is not univalent if $a>2$.
\par
Indeed, the same argument we have used shows that the function $F=H+\overline G$, where $H$ and $G$ solve the system
\[
\left\{
\begin{array}{c}
H-G= k_a\\
\omega=G^\prime/H^\prime=-l_R
\end{array}\right.\,,\quad\quad H(0)=G(0)=0\,,
\]
won't be univalent if $a>2$ either. Note that $F$ is univalent if and only if the function $f$ defined by $f(z)=-F(-z)$ is univalent. By definition, $F=K_H(1,a,-1,R)$ and using Property (iii) mentioned above, we have that $f$ equals $K_H(1,-a,1,R)=K_{-a, R}$. Therefore, we conclude that $K_{-a, R}$ is univalent if and only $0<-a\leq 2$ . This ends the proof of the theorem.
\end{pf}
\par\smallskip
In the next section of the paper we will analyze some further properties of these generalized harmonic Koebe mappings $K_{a, R}$.

\section{Extremal properties of generalized harmonic Koebe functions}\label{sec-main}

Let us again use the notation
\[
h(z)=z+\sum_{n=2}^\infty a_n(h) z^n \quad {\rm and} \quad g(z)=\sum_{n=2}^\infty b_n(g) z^n
\]
for the Taylor series expansions of functions in the canonical decomposition $f=h+\overline g$ of a sense-preserving harmonic mapping  in $\D$. We will just write $a_n(h)=a_n$ and $b_n(g)=b_n$ if there is no ambiguity.
\par
First, we would like to point out a couple of remarks.
\par
The first one is related to the coefficient $b_2$ of sense-preserving harmonic functions in the unit disk $f=h+\overline g$ normalized by $h(0)=g(0)=g'(0)=1-h'(0)=0$. Notice that for such mapping $f$, we have that $\omega= g'/h'$ is an analytic function in the unit disk that fixes the origin. Hence, we have by the Schwarz lemma that $|\omega'(0)|\leq 1$. Since we can write $g'=\omega h'$, we get
\[
|b_2|=\frac 12 |g''(0)|=\frac 12 |\omega'(0)h'(0)+\omega(0)h''(0)|=\frac 12 |\omega'(0)|\leq \frac 12.
\]
Therefore, we have that the real part of the coefficient $b_2$ is always less than or equal to $1/2$.
\par\smallskip
The second remark regarding the coefficient $a_2$ with maximum real part in the family $F_H^0$ is stated as a lemma.
\begin{lem}\label{lem-a2}
Let $F_H$ be an affine and linear invariant family of (sense-pre\-serv\-ing) harmonic mappings in $\D$. Assume that there exists $f_0=h_0+\overline{g_0}\in F_H^0$ such that
\[
\sup_{f=h+\overline g\in F_H^0} Re \{a_2(h)\}=a_2(h_0)\quad {\rm and}\quad \sup_{f=h+\overline g\in F_H^0} Re \{b_2(g)\}=b_2(g_0)\,.
\]
Then $a_2(h_0)> 1/2$.
\end{lem}
\begin{pf}
Assume first, in order to get a contradiction, that $a_2(h_0)<1/2$. (Recall that since $a_2(h_0)$ maximizes $Re\{a_2(h)\}$ in $F_H^0$, we necessarily have that $a_2(h_0)$ is a non-negative real number.) Then, we obtain that for any function $f=h+\overline g\in F_H^0$, the bound $Re\{a_2(h)\}<1/2$ holds as well.
\par
Now, let $\Phi=\Psi+\overline \Gamma$ be any function in $F_H$.
Arguing as on \cite[p. 79]{Dur-Harm}, we see that there exist $f\in F_H^0$ and $\beta\in\D$ such that $\Phi=f+\overline{\beta f}$. A straightforward computation shows that $\Psi=h+\overline\beta g$. Therefore,
\begin{eqnarray}\label{eq-lemcontrad}
\nonumber \sup_{\Phi=\Psi+\overline \Gamma\in F_H} Re\{a_2(\Psi)\} &=& \sup_{f=h+\overline g\in F_H^0, \beta\in \D} \frac{Re\{h''(0)+\overline\beta g''(0)\}}{2}\\
\nonumber &\leq& \sup_{f=h+\overline g\in F_H^0, \beta\in \D} \frac{Re\{h_0''(0)\}+ |\beta g''(0)|}{2}\\ &\leq&  a_2(h_0)+ \sup_{f=h+\overline g\in F_H^0}|b_2(g)| <1\,.
\end{eqnarray}
\par
Denote by $\mathcal H$ the family of all analytic parts of functions in $F_H$. Since $F_H$ is affine and linear invariant, it is clear that $\mathcal H$ is a linear invariant family of analytic functions. As mentioned before, Pommerenke proved that the supremum of $Re\{a_2(\Psi)\}$ for $\Psi$ in any linear invariant family of analytic functions is always bigger than or equal to $1$. This contradicts (\ref{eq-lemcontrad}) and proves that $a_2(h_0)\geq 1/2$.
\par
To finish the proof, we are to show that $a_2(h_0)=1/2$ is not possible. Again, to get a contradiction, we suppose that $a_2(h_0)=1/2$. We have two different cases.
\par
The first case is
\[
\sup_{f=h+\overline g\in F_H^0} |b_2(g)|=b_2(g_0)<\frac 12\,.
\]
Then, arguing as above, we see that the order of the linear invariant family $\mathcal H$ is strictly less than one, which is impossible.
\par
The second case to be analyzed is when
\[
 b_2(g_0)= \frac 12\,.
\]
\par
As before, denote by $\mathcal H$ the linear invariant family of analytic functions which are analytic parts of mappings in $F_H$.  The order of this family coincides with the order of the closure (with respect to the topology of locally uniform convergence) $\overline{\mathcal H}$ of $\mathcal H$. Our hypotheses show that the order $\alpha(\overline{\mathcal H})$ of $\overline{\mathcal H}$ is equal to $1$. By definition, $\alpha(F_H)=\alpha(\mathcal H)$. Hence, $\alpha(F_H)=1$.
\par
On the other hand, since $b_2(g_0)=1/2$, the dilatation $\omega=g_0'/h_0'$ of $f_0$ satisfies $\omega'(0)=2b_2(g_0)=1$. By the Schwarz lemma, we conclude that $\omega_0$ equals the identity function. In other words, we have seen that there is a harmonic mapping $f_0=h_0+\overline g_0 \in F_H^0$ with dilatation equal to the identity. Hence, following the same argument as in the proof of \cite[Thm. 3]{CHM}, we get that $\alpha(F_H)\geq 2$. This gives us the desired contradiction and ends the proof of this lemma.
\end{pf}
\par\smallskip
Note that we have just proved that if $f=h+\overline g$ maximizes simultaneously the second Taylor coefficients of both analytic and co-analytic parts of functions in $F_H^0$, then $a_2(h)-b_2(g)>0$.
\begin{thm}\label{thm-main}
Let $F_H$ be an affine and linear invariant family of (normalized) sense-preserving harmonic mappings. Assume that there is a function $f_0=h_0+\overline{g_0}\in F_H^0$ such that
\begin{equation}\label{eq-thmmain}
\sup_{f=h+\overline g\in F_H^0} Re \{a_n(h)\}=a_n(h_0)\quad  {\rm  and}\quad  \sup_{f=h+\overline g\in F_H^0} Re \{b_n(g)\}=b_n(g_0)
\end{equation}
for all $n\geq 2$. Then $f_0=K_{a, R}$, where $K_{a, R}$ is a generalized harmonic Koebe function of the form (\ref{eq-harmonicgeneralizedKoebe}). Moreover, the order of $F_H$ equals $\alpha(F_H)=a+R$.
\end{thm}
\begin{pf}
In order to simplify the notation, let us denote by $A_n=a_n(h_0)$ and $B_n=b_n(g_0)$.
\par
Assume that such an extremal function $f_0=h_0+\overline{g_0}$ exists. Using the normalizations for functions in $F_H^0$ we have $A_0=B_0=B_1=0$ and $A_1=1$, so that we can write
\[
f_0(z)=z+\sum_{n=2}^\infty A_nz^n+\overline{\sum_{n=2}^\infty B_n z^n}\,.
\]
\par
The hypotheses (\ref{eq-thmmain}) implies that the Marty relations (\ref{eq-martyan}) and (\ref{eq-martybn}) hold for all $n,\ n\geq 1$. As it was mentioned in Section~\ref{ssec-Marty}, since (\ref{eq-thmmain}) holds, all Taylor coefficients $A_n$ and $B_n$ are non-negative real numbers for $n\geq 1$. Therefore, we have that
\begin{equation}\label{eq-thmmainmartyan}
(n+1)A_{n+1} =2A_2A_n+2B_2B_n+(n-1)A_{n-1}
\end{equation}
and
\begin{equation}\label{eq-thmmainmartybn}
(n+1)B_{n+1} =2A_2B_n+2B_2A_n+(n-1)B_{n-1}
\end{equation}
occur for all such $n$\,.
\par
Now, by summing up (\ref{eq-thmmainmartyan}) and (\ref{eq-thmmainmartybn}) we get
\begin{eqnarray}\label{eq-thmmainmartyan+bn}
\nonumber  & (n+1)  \left(A_{n+1}+B_{n+1}\right)\\
& \nonumber=2(A_2+B_2)(A_n+B_n) +(n-1)(A_{n-1}+B_{n-1})
\end{eqnarray}
and hence, the function $h_0+g_0$ satisfies (\ref{eq-edo}) with $\alpha=2(A_2+B_2)$, so that $h_0+g_0$ equals the (analytic) generalized Koebe function $k_{A_2+B_2}$. (Note that the same conclusion would have been obtained by applying Proposition~\ref{prop-main}.)
\par
Now that we have proved that $h_0+g_0=k_{A_2+B_2}$, in order to determine $f_0$ completely we are to compute its dilatation $\omega$. To do so, note that
by subtracting (\ref{eq-thmmainmartyan}) and (\ref{eq-thmmainmartybn}) we get
\begin{eqnarray}\label{eq-thmmainmartyan+bn}
\nonumber  & (n+1)  \left(A_{n+1}-B_{n+1}\right)\\
& \nonumber=2(A_2-B_2)(A_n-B_n) +(n-1)(A_{n-1}-B_{n-1})
\end{eqnarray}
and hence, the function $h_0-g_0$ satisfies (\ref{eq-edo}) with $\alpha=2(A_2-B_2)$ (recall that we always have that $A_2-B_2>0$ by Lemma~\ref{lem-a2}), so that $h_0-g_0$ equals  $k_{A_2-B_2}$. Therefore,
(using also that $h_0+g_0=k_{A_2+B_2}$), we obtain
\[
h_0'(1+\omega)=k'_{A_2+B_2} \quad {\rm and}\quad h_0'(1-\omega)=k'_{A_2-B_2}\,,
\]
so that
\[
\frac{1+\omega(z)}{1-\omega(z)}=\frac{k_{A_2+B_2}'(z)}{k_{A_2-B_2}'(z)}=\left(\frac{1+z}{1-z}\right)^{2B_2}
\]
or
\[
w(z)=\frac{\left(\frac{1+z}{1-z}\right)^{2B_2}-1}{\left(\frac{1+z}{1-z}\right)^{2B_2}+1}=l_{2B_2}\,.
\]
This shows that $f_0=K_{A_2-B_2,\ 2B_2}$. (Thus, the parameters $a$ and $R$ in the statement of the theorem equal $a=A_2-B_2$ and $R=2B_2$, respectively.)
\par
To prove the second part of the theorem about the order of $F_H$, first note that $a+R=A_2+B_2$, where
\[
A_2=\max_{f=h+\overline g\in F_H^0} Re\{a_n(h)\}=\max_{f=h+\overline g\in F_H^0} |a_n(h)|
\]
and
\[
B_2=\max_{f=h+\overline g\in F_H^0} Re\{b_n(g)\}=\max_{f=h+\overline g\in F_H^0} |b_n(g)|\,.
\]
As mentioned in the proof of Lemma~\ref{lem-a2}, any function $\Phi=\Psi+\overline \Gamma \in F_H$ can be written as
\[
\Phi=f+\overline{\beta f}
\]
for some $f=h+\overline g \in F_H^0$ and $\beta\in\D$. Indeed, $\Psi=h+\overline \beta g$ , so that we get that the second Taylor coefficient of any such $\Psi$ satisfies
\[
|a_2(\Psi)|=|a_2(h+\overline \beta g)|\leq A_2+B_2\,,
\]
which proves that $\alpha(F_H)\leq A_2+B_2$. Since the function $h_0+g_0$ belongs to the closure $\overline{\mathcal H}$ of the linear family of analytic parts of functions in the family $F_H$ and, by definition, $\alpha(F_H)=\alpha(\mathcal H)(=\alpha(\overline{\mathcal H}))$, we get the desired result.
\end{pf}
\par\smallskip


\section{Some final remarks and open questions}\label{sec-final}

To finish this paper, we would like to point out several questions regarding the generalized harmonic Koebe functions defined in Section~\ref{ssec-harmkoebe}.
\par\smallskip
First of all, we want to stress that in spite of the fact that we have just considered generalized harmonic Koebe functions $K_H=K_H(\lambda, a, \mu, R)=h+\overline g$, where $h$ and $g$ solve the system
(\ref{eq-harmKoebe3}) for real parameters $\lambda$, $\mu$, and $a$ (recall that we have used the notation $K_{a,R}$ in these cases), it might be of some interest to determine what are the values of the parameters for which the corresponding generalized harmonic Koebe mapping turns out to be univalent. Though we have not been able to solve this problem in general, we can say that it is not difficult to check that by considering rotations of the form $\overline \eta K_{a,R}(\eta z)$, $|\eta|=1$, and using Theorem~\ref{thm-univalent}, one can determine different values of the parameters (that are not necessarily real) for which $K_H(\lambda, a, \mu, R)$ is one-to-one.
\par\smallskip
However, we think that it is worth mentioning that the cases studied here (this is, those functions $K_{a,R}$) seem to be important by themselves since, as was shown in Theorem~\ref{thm-main}, they are precisely the functions maximizing the Taylor coefficients of functions in affine and linear invariant families of harmonic mappings in the unit disk. Also, they maximize the growth and distortion of functions in these families. More concretely, as it is mentioned in \cite{S-S} (see also \cite[p. 97--99]{Dur-Harm}), any function $f=h+\overline g$ normalized so that $g'(0)=0$ in an affine and linear invariant family of univalent harmonic mappings with order $\alpha$ satisfies
\begin{equation}\label{eq-ultima1}
\frac{1}{2\alpha}  \left[1- \left(\frac{1-|z|}{1+|z|}\right)^\alpha\right] \leq |f(z)| \leq \frac{1}{2\alpha} \left[\left(\frac{1+|z|}{1-|z|}\right)^\alpha-1\right]\,,
\end{equation}
\begin{equation}\label{eq-ultima2}
\frac{(1-|z|)^{\alpha-1}}{(1+|z|)^{\alpha+1}} \leq |h'(z)|-|g'(z)|\,,
\end{equation}
and
\begin{equation}\label{eq-ultima3}
|h'(z)|+|g'(z)| \leq \frac{(1+|z|)^{\alpha-1}}{(1-|z|)^{\alpha+1}}\,.
\end{equation}
In principle, perhaps one would expect that functions $f$ for which equalities in (\ref{eq-ultima1}), (\ref{eq-ultima2}), or (\ref{eq-ultima3}) hold should have dilatation of the form $\omega(z)=\lambda z$ for some $|\lambda|=1$. However, there are harmonic mappings that give equality in all these inequalities having dilatation different from a rotation of the unit disk. Namely, it is not difficult to check that all functions $K_{a,R}$ with $\alpha=a+R$ (and their rotations) do it (being $R=1$ or not).
\section*{Acknowledgements}

This work is part of the first author's Ph.D. Thesis work at the Pontificia Universidad Cat\'olica de Chile under the supervision of Martin Chuaqui.
He would like to thank him for his encouragement and help.
\par
This research was developed during the first author's research stay at the University
of Eastern Finland. He also wants to use this opportunity to thank the Department of Physics and Mathematics
for its hospitality.


\end{document}